\documentclass[final,seceqn]{elsarticle}

\usepackage{mathrsfs,amsfonts,amssymb,latexsym,amsmath,amscd}
\usepackage{graphicx,epsf,amsbsy,bm}

\journal{Applied Mathematics and Computation}

\usepackage{hyperref}
\hypersetup{
colorlinks=true,
linkcolor=blue,
filecolor=magenta,
urlcolor=cyan,
}
\newtheorem{theorem}{Theorem}

\newtheorem{corollary}{Corollary}

\newenvironment{proof}[1][Proof.]{\begin{trivlist}
\item[\hskip \labelsep {\bfseries #1}]}{\end{trivlist}}

\newcommand{\be}{\begin{equation}}
\newcommand{\ee}{\end{equation}}
\newcommand{\bea}{\begin{eqnarray}}
\newcommand{\eea}{\end{eqnarray}}
\newcommand {\mm}[1]{\quad\mbox{#1}\quad}

\newcommand{\R}{{\Re}}
\newcommand {\lr}[1]{\left({#1}\right)}

\newcommand {\lrs}[1]{\left[{#1}\right]}

\newcommand {\eq}[1]{(\ref{#1})}
\newcommand {\TF}[6]
{{}_3F_2\!\left[\!\!\begin{array}{c}
#1,\;#2,\;#3;\\
\\
#4,\;#5;
\end{array}
\!\;#6\!\right]}
\newcommand {\FG}[4]
{{}_2F_1\!\left[\!\!\begin{array}{c}
#1,\;#2;\\
\\#3;
\end{array}
\!\;#4\!\right]}

\begin{document}
\begin{frontmatter}
\title{The Clausenian hypergeometric function $_3F_2$ with unit argument
and negative integral parameter differences}

\author[mymainaddress,mysecondaryaddress]{M. A. Shpot\corref{cor1}\fnref{myfootnote}}
\cortext[cor1]{Corresponding author}
\fntext[myfootnote]{\url{http://researchgate.net/profile/Mykola_Shpot}}
\address[mymainaddress]
{Institute for Condensed Matter Physics, 79011 Lviv, Ukraine}
\address[mysecondaryaddress]
{Fakult\"at f\"ur Physik, Universit\"at Duisburg-Essen,
D-47048 Duisburg, Germany}

\author[HMSaddress]{H. M. Srivastava\fnref{HMSfootnote}}
\fntext[HMSfootnote]{\url{http://www.math.uvic.ca/faculty/harimsri/}}
\address[HMSaddress]
{Department of Mathematics and Statistics, University of Victoria,
Victoria, British Columbia V8W 3R4, Canada}

\begin{abstract}
New explicit as well as manifestly symmetric three-term summation
formulas are derived for the Clausenian hypergeometric series
$_3F_2(1)$ with negative integral parameter differences.
Our results generalize and naturally extend several similar relations
published, in recent years, by many authors.
An appropriate and useful connection is established with the quite
underestimated 1974 paper by P. W. Karlsson.
\end{abstract}

\begin{keyword}
\texttt{Generalized hypergeometric functions}
\sep\texttt{Clausen series}
\sep\texttt{summation formulas}
\sep\texttt{hypergeometric identities}
\sep\texttt{integral parameter differences}
\MSC[2010] 33C20\sep 33C05\sep 33C90\sep 33B15\sep 33B20

\end{keyword}
\end{frontmatter}

\section{Introduction}\label{sec1}

The main object of the present paper is to derive two elegant and
manifestly symmetric summation formulas \eq{ZY} and \eq{ZX} for
Clausen's series $_3F_2$ with unit argument (see, e.g.,
\cite{Bailey,9,Rainv,Slater,Luke,SriMan,SK})
\be\label{VB}
\TF{a}{b}{c}{b+1+m}{c+1+n}{1} \qquad (m,n\in \mathbb{N}_0:=\mathbb{N}\cup\{0\}),
\ee
and \emph{negative} integral parameter differences given by the arbitrary non-negative
integers $m$ and $n$, $\mathbb{N}$ being (as usual) the set of {\it positive} integers.
We also aim at supporting interest in generalized hypergeometric functions $_pF_q (z)$
of $p$ numerator and $q$ denominator parameters (and with argument $z$)
especially of the type exemplified above.

Our summation formulas \eq{ZY} and \eq{ZX} match very
well the recent trend of finding new relationships for
generalized hypergeometric functions.
In fact, they are immediate and natural generalizations of more special
formulas suggested, a decade ago,
by Milgram \cite{Milg04,Milg10,Milg447,Milg07c},
which were further proved and employed by Miller and Paris
\cite{MP12} and Rathie and Paris \cite{RP13p} quite recently.

A substantial and wide-spread progress has been recently achieved in the
classical field of investigating the generalized hypergeometric
functions $_pF_q(z)$ and producing various relationships between them.
Very often the studied functions contain, in different ways,
integers in their numerator and denominator parameters
(see, for example, \cite{LGR96,Lew97,RBK05,KratRiv06,RR08,MS10,MP11,Chu12,RP13,KRP14T}).

More studied, however, are functions with \emph{positive}
integral parameter differences
in pairs of their numerator and denominator parameters
\cite{MS10,MP11,MP13,RP13,KRP14T} just in the spirit
of very well-known early papers by Karlsson \cite{K71}
and Minton \cite{Mint70}. It may be of interest to recall here
the following known reduction formula (see, for example, \cite{K71}, \cite{Sri1-1},
and also \cite[p. 1080]{RSri} and the references to more general results on hypergeometric
reduction formulas, which are cited in \cite{RSri}):
\begin{align}\label{Kar}
& \;_pF_q \left[ {{\begin{array}{c}
 {b_1+m_1,\cdots,b_r+m_r, a_{r+1},\cdots,a_p;}  \\
\\{b_1,\cdots,b_r,b_{r+1},\cdots,b_q;}
\end{array}}z} \right]\notag \\[3mm]
& \qquad \qquad =\sum_{j_1=0}^{m_1}\cdots\sum_{j_r=0}^{m_r}\Lambda(j_1,\cdots,j_r)
z^{J_r}\;_{p-r}F_{q-r} \left[ {{\begin{array}{rr}
 {a_{r+1}+J_r,\cdots,a_p+J_r;}  \\
\\{b_{r+1}+J_r,\cdots,b_q+J_r;}
\end{array}}z} \right]
\end{align}
$$(r\leqq \min\{p,q\};\; p,q\in \mathbb{N}_0;\; p<q+1 \quad \text{when}
\quad z\in \mathbb{C};\;
p=q+1 \quad \text{when} \quad |z| < 1),$$
where
$$J_r := j_1+\cdots + j_r$$
and
$$\Lambda(j_1,\cdots,j_r)=\binom{m_1}{j_1}\cdots\binom{m_r}{j_r}
\frac{(b_2+m_2)_{J_1}\cdots (b_r+m_r)_{J_{r-1}}\;(a_{r+1})_{J_r}\cdots (a_p)_{J_r}}
{(b_1)_{J_1}\cdots (b_r)_{J_r}(b_{r+1})_{J_r}\cdots (b_q)_{J_r}}.$$

The general hypergeometric identity \eqref{Kar} was proved by
Karlsson \cite{K71} and (in two markedly different {\it simpler} ways) by
Srivastava \cite{Sri1-1}. More interestingly, various generalizations and {\it basic}
(or $q$-) extensions of the hypergeometric identity \eqref{Kar} can be found in
several sequels to the works by Karlsson \cite{K71} and Srivastava \cite{Sri1-1}
(see, for example, \cite{Panda1}). Reference \cite{Panda2}, on the other hand, contains
{\it further} general results stemming from the hypergeometric identity \eqref{Kar} including
{\it multivariable} generalizations. Furthermore, Karlsson's proof of the Karlsson-Minton
summation formula (see, for details, \cite{K71}; see also \cite{Mint70}, \cite{MS10} and
\cite[p. 1080, Equation (20)]{RSri}) was based upon the hypergeometric reduction formula \eqref{Kar}.

There is another obscure and seemingly forgotten paper by
Karlsson \cite{K74} in which similar results have been obtained for
generalized hypergeometric functions $_pF_q(z)$ with \emph{negative}
integral parameter differences.
The papers \cite{Milg04,Milg07c} and \cite{MP12}
mentioned at the beginning, as well as the present one,
discuss the summation formulas for the functions
$_3F_2(1)$ that belong to the same category.

Our motivation in doing this work stems from calculations
\cite{DS00,SD01,SPD05,RDS11,SP12}
in the field theory of Lifshitz points \cite{HLS75},
where $_3F_2$ functions of the type indicated in \eq{VB}
appear as a part of the expansion coefficients of certain
important functions
(see, for explicit formulas, \cite[Eqs. (5.69), (5.71)]{RDS11}).
Such expansions appear as a result of a term-by term integration
of special \cite{Sh07} Appell functions
\cite{Appell26,Bailey,Slater,SriMan,SK}.
Owing to global universal features of mathematical description
of the underlying systems with anisotropic scaling,
similar expansions are expected to inevitably appear in a very
broad class of statistical physics and (Lorentz violating)
high energy theories as discussed in a review section of
\cite{SP12}. On the other side, a review and further references can be found
in \cite{Sot11}.
Owing to numerous potential applications, both in theoretical physics and
mathematics \cite{MP12,RP13p}, we believe that functions
\eq{VB} or the related ones, deserve to be studied in a best way.

In the following section, we explicitly write down the previous results of \cite{K74},
\cite{Milg04,Milg07c} and \cite{MP12}, which will be
needed for establishing the necessary contacts and connections with the present work.

\section{Background results}\label{BR}

In 1974, Karlsson \cite{K74} derived a quite general reduction formula for
generalized hypergeometric functions
$_pF_q(z)$ with generic negative integral parameter differences
and for $p\leqq q+1$. In the case when $p=q+1=3$, his equation (6) in \cite{K74}
may be written as follows:\footnote{The restriction $\mu\ne 1$ appearing in (7)
has to be read as $\mu\ne i$.}
\begin{align}\label{KA}
&\TF{a}{b}{c}{b+1+m}{c+1+n}{z}\frac{m!\,n!}{(b)_{m+1}(c)_{n+1}}=
\sum\limits_{i=0}^m\sum\limits_{j=0}^n\frac{(-m)_i(-n)_j}{i!\,j!}\nonumber\\
&\qquad\cdot
\left({\frac{1}{(c+j)(b-c+i-j)}\;\FG{a}{c+j}{c+j+1}{z}+
\frac{1}{(b+i)(c-b+j-i)}\;\FG{a}{b+i}{b+i+1}{z}}\right),
\end{align}
provided that no denominator parameter equals zero or a negative integer
and $|\arg(1-z)|<\pi$.
Here, and throughout this paper, $m,n\in\mathbb N_0$ are arbitrary
non-negative integers. The Pochhammer symbol $(\lambda)_n$ is given by
(see, e.g., \cite[Ch. 1]{SriMan} and \cite{MS10})
\be\label{PL}
(0)_0:=1 \qquad \text{and} \qquad (\lambda)_n \equiv\lambda(\lambda+1) \cdots (\lambda+n-1)\qquad
(\lambda \in \mathbb C;\; n\in\mathbb{N})
\ee
and, in general, by
\be\label{PG}
(\lambda)_n =
\frac{\Gamma(\lambda+n)}{\Gamma(\lambda)}\qquad (n\in\mathbb{Z}).
\ee
Finally, $_2F_1(a,b;c;z)$ is a Gauss hypergeometric function and
$\Gamma(z)$ is the usual Euler Gamma function
(see, e.g., \cite{Bailey,Slater,9,SriMan}).

At unit argument $z=1$, the Gauss hypergeometric functions on the right-hand side of \eqref{KA}
are summed by applying the celebrated Gauss summation theorem \cite[Sec. 1.3]{Bailey}:
\begin{equation}\label{GS}
\FG{a}{b}{c}{1}=\frac{\Gamma(c)\Gamma(c-a-b)}{\Gamma(c-a)\Gamma(c-b)}
\qquad \big(\R(c-a-b)>0\big).
\end{equation}
We thus find from \eq{KA} that
\be\label{KB}
\TF{a}{b}{c}{b+1+m}{c+1+n}{1}\frac{m!\,n!}{(b)_{m+1}(c)_{n+1}}=
\sum\limits_{i=0}^m\sum\limits_{j=0}^n\frac{(-m)_i(-n)_j}{i!\,j!}
\left({\frac{B(1-a,c+j)}{b-c+i-j}+\frac{B(1-a,b+i)}{c-b+j-i}}\right).
\ee
Here
\be\label{BE}
B(a,b)=\frac{\Gamma(a)\Gamma(b)}{\Gamma(a+b)}\qquad (a,b\ne 0,-1,-2,\cdots)
\ee
denotes the familiar Beta function.

In 2004, Milgram \cite[Eq. (11)]{Milg04} suggested the following summation formula
for Clausen's $_3F_2(1)$ series:
\begin{align}\label{OM}
&\TF{a}{b}{c}{b+n}{c+1}{1}=\frac{(b)_n\Gamma(c+1)\Gamma(1-a)}
{(b-c)_n\Gamma(c+1-a)}\notag \\
& \qquad\quad +c\Gamma(b+n)\Gamma(c-b+1-n)
\sum_{\ell=0}^{n-1}\frac{\Gamma(n-\ell-a)(-1)^\ell}
{\Gamma(b+n-a-\ell)\Gamma(n-\ell)\Gamma(c-b-n+2+\ell)}.
\end{align}
He further reproduced the summation formula \eqref{OM} in slightly different forms in
\cite{Milg07c,Milg10,Milg447}.

Quite recently, Miller and Paris \cite{MP12} re-derived the summation
formula \eq{OM} in further two equivalent forms.
Their equation (1.6) reads as follows:
\be\label{MP}
\TF{a}{b}{c}{b+n}{c+1}{1}=\frac{c\Gamma(1-a)(b)_n}{(b-c)_n}
\left(\frac{\Gamma(c)}{\Gamma(1+c-a)}-\frac{\Gamma(b)}{\Gamma(1+b-a)}
\sum_{k=0}^{n-1}\frac{(1-a)_k(b-c)_k}{(1+b-a)_k\,k!}\right),
\ee
where it is expressed in terms of the partial sum of a Gauss
hypergeometric function $_2F_1$ of unit argument,
$$
\FG{a}{b}{c}{1}_n\equiv\sum_{k=0}^{n}\frac{(a)_k(b)_k}{(c)_k\,k!}.
$$
By Lemma 2 of \cite{MP12}, which reads
\be\label{MT}
\FG{a}{b}{c}{1}_n=\frac{(1+b)_n}{n!}\;\TF{-n}{b}{c-a}{1+b}{c}{1},
\ee
the equation \eq{MP} can be equivalently written in
terms of the terminating Clausen's series $_3F_2(-n,b,c;d,e;1)$
(see, for details, \cite[Eq. (3.1)]{MP12}).

In order to facilitate the comparison with our further results, we
write down two variants of \eq{MP} for $n\mapsto n+1$.
Following the choice $a\mapsto 1-a$ and $b\mapsto b-c$ with $c\mapsto 1+b-a$
of \cite{MP12} in applying \eq{MT} to \eq{MP}, we obtain
\be\label{FA}
\TF{a}{b}{c}{b+1+n}{c+1}{1}=c\,(b)_{n+1}
\left(\frac{B(1-a,c)}{(b-c)_{n+1}}-
\frac{B(1-a,b)}{n!}\,\TF{-n}{b}{b-c}{1+b-a}{1+b-c}{1}\right),
\ee
which exactly matches \cite[Eq. (3.1)]{MP12} in view of \eq{BE}.
Alternatively, upon setting $a\mapsto b-c$ and $b\mapsto 1-a$ with the same $c$
as above (that is, with $c\mapsto 1+b-a$), we arrive at the following result:
\begin{align}\label{FB}
&\TF{a}{b}{c}{b+1+n}{c+1}{1}=\frac{c(b)_{n+1}}{(b-c)_{n+1}}\notag \\
&\qquad \quad \cdot \left(B(1-a,c)-B(1-a,b)\,\frac{(2-a)_n}{n!}\,
\TF{-n}{1-a}{1-a+c}{2-a}{1-a+b}{1}\right).
\end{align}

In closing this section, let us only note that the $_3F_2$ functions
in the last two equations \eqref{FA} and \eqref{FB} are related by means
of the following two-term Thomae
transformation (see, e.g., \cite[Entry (7.4.4.1)]{PBM3}):
\be\label{TH}
\TF{a}{b}{c}{d}{e}{1}=
\frac{\Gamma(d)\,\Gamma(p-c)}{\Gamma(p)\,\Gamma(d-c)}
\;\TF{e-a}{e-b}{c}{p}{e}{1}\qquad (p:=d+e-a-b),
\ee
which leaves the last pair of the numerator and denominator parameters
(that is, $c$ and $e$) unaltered.

\section{Summation theorems}\label{sec4}

In this section, we prove two theorems that generalize the results quoted
in the Section \ref{BR}. Theorem 1 gives a
symmetric variant of \eq{FA} and \eq{FB} in which the second denominator
parameter $c+1$ is replaced by $c+1+n$.
Our Theorem \ref{E2} will go a step further by allowing the
negative integral differences in two pairs of parameters
of $_3F_2$ series to be independent.

\begin{theorem}\label{E1}
For arbitrary non-negative integers $n\in\mathbb N_0$
and complex numbers $a,\ b,\ c\in\mathbb C,$
\begin{align}\label{TT}
&\TF{a}{b}{c}{b+1+n}{c+1+n}{1}\frac{(c-b)_{n+1}\,n!}{(b)_{n+1}(c)_{n+1}}\notag \\
&\qquad\quad
= B(1-a,b)\,\frac{(1-a)_n}{(1+b-a)_n}\,\TF{-n}{b}{1+n}{1+b-c}{a-n}{1}\notag \\
&\qquad\qquad\qquad\quad
+(-1)^{1+n}\,B(1-a,c)\,\frac{(2-a+n)_n}{(1-a+c)_n}\;
\TF{-n}{1-a+b+n}{1+n}{1+b-c}{2-a+n}{1},
\end{align}
provided that $\R(2-a+2n)>0$.
\end{theorem}

\begin{proof}
In transforming the left-hand side of the assertion \eqref{TT},
we use the following Thomae three-term relation for $_3F_2$
at unit argument quoted by Bailey \cite[p. 21, Eq. (1)]{Bailey},
which can be also found in \cite[Entry (7.4.4.4)]{PBM3}:
\footnote{In this reference \cite{PBM3}, this is the \emph{next} entry to that
employed by Miller and Paris \cite{MP12} in their proof of the summation
formula \eq{MP}.}
\begin{eqnarray}\label{Thomae2}
&&\TF{a}{b}{c}{e}{f}{1}=
\frac{\Gamma(e)\Gamma(e-a-b)}{\Gamma(e-a)\Gamma(e-b)}\;
\TF{a}{b}{f-c}{a+b-e+1}{f}{1}\nonumber
\\
&&\qquad\qquad
+\frac{\Gamma(e)\Gamma(f)\Gamma(a+b-e)\Gamma(e+f-a-b-c)}
{\Gamma(a)\Gamma(b)\Gamma(f-c)\Gamma(e+f-a-b)}
\;\TF{e-a}{e-b}{e+f-a-b-c}{e-a-b+1}{e+f-a-b}{1},
\end{eqnarray}
where we have $\R(e+f-a-b-c)>0$ for convergence of the $_3F_2(1)$
series on the left-hand side
and $\R(1+c-e)>0$ for convergence of both of the $_3F_2(1)$ series on the right-hand side.
Up to their coefficients, the two resulting $_3F_2(1)$ functions are given by
\be\label{NF}
\TF{a}{b}{1+n}{a-n}{1+c+n}{1}\qquad\mm{and}\qquad
\TF{1-a+b+n}{1+n}{2-a+2n}{2-a+n}{2-a+c+2n}{1},
\ee
which converge when $\R(c-b)>n$.

The crucial property of the applied transformation \eqref{Thomae2}
is that the first of the functions in \eq{NF}
has a pair of the upper and lower parameters $(a,a-n)$
and the second one has $(2-a+2n,2-a+n)$.
In both cases the denominator parameter differs from the numerator parameter
by a negative integer $-n$.
Thus, for any non-negative finite integer $n<\R(c-b)$, each of the
new $_3F_2$ functions is reducible to a finite sum of products
of Euler Gamma functions in light of the following relation \cite[Entry (7.4.1.2)]{PBM3}:
\be\label{PF}
\TF{a}{b}{c}{a-n}{d}{z}=\frac{1}{(1-a)_n}\sum_{p=0}^n (-z)^p
\binom{n}{p}(1-a)_{n-p}\frac{(b)_p(c)_p}{(d)_p}\,\FG{b+p}{c+p}{d+p}{z}
\ee
with $z=1$ as well as the Gauss summation theorem \eq{GS}.

Expressing the binomial coefficients in \eq{PF} via
(see, e.g., \cite[p. 22, Eq. (16)]{SriMan})
$$
\binom{n}{p}=\frac{(-1)^p(-n)_p}{p!}
$$
and after some algebra, we obtain a linear combination of finite sums as follows:
$$
\sum_{p=0}^n\frac{(-n)_p(1+n)_p}{(a-n)_p\,p!}\frac{(b)_p}{(1+b-c)_p} \qquad
\mm{and} \qquad
\sum_{p=0}^n\frac{(-n)_p(1+n)_p(1-a+b+n)_p}
{(2-a+n)_p\,p!}\frac{1}{(1+b-c)_p}.
$$
Identifying each of these finite sums with terminating Clausen's series $_3F_2(1)$,
we derive
\bea\label{TO}
&&\TF{a}{b}{c}{b+1+n}{c+1+n}{1}\frac{(c-b)_{n+1}\,n!}{(b)_{n+1}(c)_{n+1}}\nonumber
\\&&\nonumber\qquad\quad
=\Gamma(b)\;\frac{\Gamma(1-a+n)}{\Gamma(1-a+b+n)}
\TF{-n}{b}{1+n}{1+b-c}{a-n}{1}
\\&&\qquad\qquad\quad
+\Gamma(c)\;
\frac{\Gamma(-1+a-n)\Gamma(2-a+2 n)}{\Gamma(a)\Gamma(1-a+c+n)}\,
\TF{-n}{1-a+b+n}{1+n}{1+b-c}{2-a+n}{1}.
\eea

Finally, upon using the reflection formula \eq{RF} to transform
$\Gamma(-1+a-n)$ in the second term, if we rearrange
the involved Gamma functions, we obtain
the summation formula recorded in \eq{TT}.
\end{proof}

\noindent
{\bf Remark 1.} By changing the numerator parameter
$a\mapsto a+2n$ in \eq{TO}, we are led easily to Corollary \ref{cor1} below.

\begin{corollary} \label{cor1}
The following summation formula$:$
\bea\label{TC}
&&\TF{a+2n}{b}{c}{b+1+n}{c+1+n}{1}\frac{(c-b)_{n+1}\,n!}{(b)_{n+1}(c)_{n+1}}\nonumber
\\&&\qquad\quad =\Gamma(b)\;
\frac{\Gamma(1-a-n)}{\Gamma(1-a+b-n)}\;
\TF{-n}{b}{1+n}{1+b-c}{a+n}{1}\nonumber
\\&&\qquad\qquad\quad +\Gamma(c)\;
\frac{\Gamma(2-a)\Gamma(a+n-1)}{\Gamma(a+2n)\Gamma(1-a+c-n)}
\;\TF{-n}{1-a+b-n}{1+n}{1+b-c}{2-a-n}{1}
\eea
holds true when $\R(2-a)>0$.
\end{corollary}

As already mentioned in Section \ref{sec1}, there are physical
applications (see, for example, \cite[Eqs. (5.69)--(5.73)]{RDS11} and \cite{SP12,SP14}),
where reduction relations of this form are relevant and potentially useful.

\vskip 2mm

The following Theorem \ref{E2} extends the result \eq{TT} of Theorem \ref{E1}.

\begin{theorem}\label{E2}
For arbitrary non-negative integers
$m\in\mathbb N_0$ and $n\in\mathbb N_0,$
and for complex parameters $a,\ b,\ c\in\mathbb C,$

\be\label{MN}
\TF{a}{b}{c}{b+1+m}{c+1+n}{1}
\frac{(c-b)_{n+1}}{(b)_{m+1}(c)_{n+1}}=T^{(1)}_{m,n}+T^{(2)}_{m,n} \qquad \big(\R(2-a+m+n)>0\big),
\ee
where
\be\label{T1}
T^{(1)}_{m,n}=B(1-a,b)\,\frac{(1-a)_m}{(1+b-a)_m\,m!}\;
\TF{-m}{b}{1+n}{1+b-c}{a-m}{1}
\ee
and
\be\label{T2}
T^{(2)}_{m,n}=(-1)^{1+m}\,B(1-a,c)\,\frac{(2-a+m)_n}{(1-a+c)_n\,n!}\;
(c-b)_{n-m}\;\TF{-n}{1-a+b+m}{1+m}{2-a+m}{1+b-c+m-n}{1},
\ee
provided that $\R(2-a+m+n)>0$.
\end{theorem}
The proof of Theorem \ref{E2} proceeds precisely along the same lines as those of Theorem \ref{E1}.

\vskip 2mm

It is fairly straightforward to see that, in its special case when $m=n$, the summation formula
\eq{MN} reduces rather trivially to \eq{TT}.

\section{Special cases and consequences}

For $n=0$, the equation \eq{TT} coincides with that of Milgram \eq{OM}
and Miller and Paris \eq{MP} in the special case when
$n=1$ and yields the following known result
\cite[Entry (7.4.4.16)]{PBM3} with
$a$ and $c$ interchanged:
$$
\TF{a}{b}{c}{b+1}{c+1}{1}=\frac{bc}{c-b}\;\Gamma(1-a)
\lrs{\frac{\Gamma(b)}{\Gamma(1-a+b)}-\frac{\Gamma(c)}{\Gamma(1-a+c)}}.
$$
Another known result \cite[Entry (7.4.4.17)]{PBM3} with
$a\leftrightarrow c$ is the special case of \eq{TT} when $n=1$.
The summation formula \eq{TT} gives access to simple generalizations of
\cite[Entries (7.4.4.16) and (7.4.4.17)]{PBM3}
with arbitrary equal integral enhancements of the
denominator parameters.

As $n\geqq 1$, the formula \eq{TT}
cannot directly match the equations \eq{FA} and \eq{FB}
for symmetry reasons: In this case, the parameter differences
in $_3F_2(1)$ become asymmetric therein.
We next show that these formulas, and hence \eq{OM} and \eq{MP},
follow from the result \eq{MN} both for $m=0$ and $n=0$.

Let us first put $m=0$ in \eq{MN} to \eq{T2}.
In this case, we have to deal with the Clausenian hypergeometric function
$$
\TF{a}{b}{c}{b+1}{c+1+n}{1}
$$
on the left-hand side. In order to compare it with \eq{OM} and \eq{MP} or with
\eq{FA} and \eq{FB}, the parameters $b$ and $c$ have to be interchanged.

At $m=0$, the $_3F_2$ function in $T^{(1)}_{m,n}$ given by \eqref{T1} reduces to $1$,
and its factor, together with the one of \eq{MN}, trivially combine to
the first term on the right-hand side in \eq{FA} or \eq{FB} with
$b\leftrightarrow c$.
The remaining non-trivial $_3F_2$ function from $T^{(2)}_{m,n}$ given by \eqref{T2} simplifies to
the following form:
\be\label{QW}
\TF{-n}{1-a+b}{1}{1+b-c-n}{2-a}{1}.
\ee
It does not match directly any of the $_3F_2$ functions from
\eq{FA} or \eq{FB} because of a subtraction $-n$
in one of the denominator parameters. To get rid of this term, we
use the relation \cite[Entry (7.4.4.85)]{PBM3}. When read in the reverse
direction, it can be written as follows:
\be\label{RD}
\TF{-n}{a}{b}{c-n}{d}{1}=\frac{(1+a-c)_n}{(1-c)_n}\;
\TF{-n}{a}{d-b}{1+a-c}{d}{1}.
\ee
The evident choice $a\mapsto1-a+b$ and $b\mapsto1$ transforms the $_3F_2$ function
from \eq{QW} via
\be\label{QE}
\TF{-n}{1-a+b}{1}{1+b-c-n}{2-a}{1}=\frac{(1-a+c)_n}{(c-b)_n}\;
\TF{-n}{1-a}{1-a+b}{2-a}{1-a+c}{1}.
\ee
With $b\leftrightarrow c$, the last hypergeometric function in \eqref{QE}
matches that in the second term of \eq{FB}.
Collecting its factors from \eq{QE}, \eq{T2}, and \eq{MN}, it is a
matter of simple algebra to bring them to the form required by \eq{FB}.

If we now take $n=0$ in the Eqs. \eq{MN}--\eq{T2}, the connection
with \eq{FA} is established with the same set of parameters
$a,\ b$, and $c$. This time we have to care
only about a simple exchange of $m$ and $n$.

Again, at $n=0$, the $_3F_2$ function from $T^{(1)}_{m,n}$ given by \eqref{T1} transforms to a
``correct'' one from \eq{FB} on applying \eq{RD} with
$a\mapsto1,\ b\mapsto b,\ c\mapsto a,$ and $d\mapsto 1+b-c$,
while the one from $T^{(2)}_{m,n}$ given by \eqref{T2} reduces to $1$. Matching the factor
at the transformed function is again simple.
Reducing $T^{(2)}_{m,0}$ to the first term in \eq{FB}
requires the use of the following transformation formula
(see, e.g., \cite[p. 22, Eq. (19)]{SriMan}):
$$
(\lambda)_{-m}=\frac{(-1)^m}{(1-\lambda)_m}, \qquad
n\in\mathbb N,\qquad \lambda\notin\mathbb Z
$$
for Pochhammer symbols in \eq{PG} with $\lambda\mapsto c-b$.

\vskip 2mm

\noindent
{\bf Remark 2.}
We note in passing that the relation \cite[Entry (7.4.4.86)]{PBM3}
rewritten similarly as \eq{RD},
\be\label{RE}
\TF{-n}{a}{b-n}{c-n}{d-n}{1}=\frac{(1+a-c)_n(1-b)_n}
{(1-c)_n(1-d)_n}\;\TF{-n}{d-b}{1-c}{1+a-c}{1-b}{1},
\ee
may be also useful in practical calculations.

For instance, the function $G_{m,k}(t)$ appearing in \cite[Eq. (6)]{KRP14T}
can be reduced with the help of \eq{RE} via
$$
\TF{-m+k}{t+k}{c-a-b-m}{c-a-m+k}{c-b-m+k}{1}\quad \longrightarrow\quad
\TF{-m+k}{1+a-c}{c-b-t-m}{1-b-m}{1-t-m}{1}
$$
to a form where the summation index $k$ appears only once in
a numerator parameter.

\section{Results in the final form}

In the preceding sections, we learned that the standard identities
\eq{TH} and \eq{RD} are useful in dealing with the involved functions.
With this in mind, we use first \eq{RD} to transform the
$_3F_2$ functions appearing in \eq{T1} and \eq{T2}.
For the first of them, with the choice $a\mapsto b$ and $b\mapsto1+n$, we obtain
\be\label{QV}
\TF{-m}{b}{1+n}{1+b-c}{a-m}{1}=\frac{(1-a+b)_m}{(1-a)_m}\;
\TF{-m}{b}{b-c-n}{1+b-a}{1+b-c}{1}.
\ee
Similarly, for the second one, we have
\be\label{QD}
\TF{-n}{1-a+b+m}{1+m}{2-a+m}{1+b-c+m-n}{1}=\frac{(1-a+c)_n}{(c-b-m)_n}\;
\TF{-n}{1-a}{1-a+b+m}{1+c-a}{2-a+m}{1}.
\ee
The $_3F_2$ function in \eq{QV} looks very good: It is a natural
generalization of that in \eq{FA}.
In our hope to obtain a more symmetric expression for the function appearing in
\eq{QD}, we transform it via \eq{TH}:
\be\label{QQ}
\TF{-n}{1-a}{1-a+b+m}{1+c-a}{2-a+m}{1}=
\frac{(1+c-b)_n}{(2-a+m)_n}\;\TF{-n}{c}{c-b-m}{1+c-a}{1+c-b}{1}.
\ee
Now the $_3F_2$ functions in \eq{QD} and \eq{QQ} are symmetric
with respect to the interchange
$b\leftrightarrow c$ and $m\leftrightarrow n$,
which is quite satisfactory. Using the last three equations
in \eq{MN}--\eq{T2} we obtain the {\color{red} final result}:

\begin{align}\label{ZY}
&\TF{a}{b}{c}{b+1+m}{c+1+n}{1}\frac1{(b)_{m+1}(c)_{n+1}}\nonumber\\
&\qquad\quad
=\frac{B(1-a,b)}{(c-b)_{n+1}\,m!}\;\TF{-m}{b}{b-c-n}{1+b-a}{1+b-c}{1}
+\frac{B(1-a,c)}{(b-c)_{m+1}\,n!}\;\TF{-n}{c}{c-b-m}{1+c-a}{1+c-b}{1}.
\end{align}

Upon setting $m=n$ in \eqref{ZY}, we are led to the following

\begin{corollary}
The summation formula $\eq{ZX}$ follows from $\eqref{ZY}$ for $m\mapsto n$:
\end{corollary}

\begin{align}\label{ZX}
&\TF{a}{b}{c}{b+1+n}{c+1+n}{1}\frac{n!}{(b)_{n+1}(c)_{n+1}} \notag \\
&\qquad\quad
=\frac{B(1-a,b)}{(c-b)_{n+1}}\;\TF{-n}{b}{b-c-n}{1+b-a}{1+b-c}{1}
+\frac{B(1-a,c)}{(b-c)_{n+1}}\;\TF{-n}{c}{c-b-n}{1+c-a}{1+c-b}{1}.
\end{align}

The Clausenian series on the left-hand side of \eq{ZY} and \eq{ZX}
converge when $$\R(2-a+m+n)>0 \qquad \text{and} \qquad \R(2-a+2n)>0,$$
respectively.
This implies that both of the summation formulas make sense
for generic $m\in\mathbb N_0$ and $n\in\mathbb N_0,$ provided that $\R(a)<2$.
Under these conditions, all functions at both sides are well defined,
and, by the principle of analytic continuation, the restriction
$n<\R(c-b)$ imposed in intermediate calculations (see \eq{NF})
can be removed.

\vskip 3mm

\noindent
{\bf Remark 3.}
Additional restrictions (for example, $a\ne 1$)
imposed by \eq{BE} on the Beta functions
in the right-hand sides of \eq{ZY} and \eq{ZX}
do not shrink the applicability region of both equations.
Singularities, which arise when the arguments of these
Beta functions approach negative integers mutually cancel in
the whole two-term combinations.
Each of such dangerous cases has to be treated separately,
in a manner similar to that of \cite{RBK05,Milg07c}:
In these references, the special case when $m=0$ and $b\to c$ has been
considered for arbitrary non-negative integer $n$.
It is not complicated to see that, when $a=1$,
the $_3F_2$ functions on the right-hand sides of \eq{ZY} and \eq{ZX}
reduce to the Gauss functions $_2F_1$, which can be summed via the
Gauss summation theorem
\eq{GS}, and the whole resulting combinations at singular
$B(1-a)$ vanish as they should.
Proceeding similarly as in \cite{RBK05,Milg07c} would lead us to a finite
limit as $a\to 1$ for the relations \eq{ZY} and \eq{ZX}.

\vskip 3mm

Let us proceed with a proposition, which
gives a simple demonstration of Karlsson's result \eq{KB},
followed by several remarks.

\noindent
{\bf Proposition 1.}
{\it The double-sum representation of $\eq{KB}$ can be expressed as a three-term
summation formula $\eq{ZY}$ for the Clausenian hypergeometric function $_3F_2$.}

\begin{proof}
In both of Karlsson's formulas \eq{KA} and \eq{KB},
the summations over $i$ in the first terms and over $j$ in the second terms
are the same up to notations and they are given by
\be
\sum_{k=0}^m\frac{(-m)_k}{(a+k)\,k!}=
\frac1{a}\sum_{k=0}^m\frac{(-m)_k(a)_k}{(a+1)_k\,k!}=
\frac1{a}\,\FG{-m}{a}{a+1}{1}=
\frac{\Gamma(a)\,m!}{\Gamma(a+1+m)},
\ee
where we have used the following simple property of the Pochhammer symbol in \eq{PL}:
$$(a)_{n+1}=a(a+1)_n=(a)_n(a+n)
$$
as well as the Gauss summation theorem \eq{GS}. We take here $a=b-c-j$
and use the familiar reflection formula (see, for example,
\cite[Ch. 2.17]{Rainv} and \cite[Ch. 1]{SriMan})
\be\label{RF}
\Gamma(a-j)\Gamma(1-a+j)=\frac\pi{\sin\pi(a-j)}=
\frac{\pi\,(-1)^j}{\sin\pi a}=(-1)^j\,\Gamma(a)\Gamma(1-a)
\qquad (a\notin \mathbb{Z}).
\ee
It gives an analytic continuation for the Euler Gamma function and implies that
$$
\Gamma(a-j)=
\frac{(-1)^j\,\Gamma(a)}{(1-a)_j}\qquad  (a\notin \mathbb{Z}).
$$
Hence we obtain
\be\label{SU}
\sum_{k=0}^m\frac{(-m)_k}{(b-c-j+k)\,k!}=
\frac{m!}{(b-c)_{m+1}}\;\frac{(c-b-m)_j}{(1+c-b)_j}
\ee
for the inner sums in the first terms of \eq{KA} and \eq{KB}.
Moreover, just as we mentioned above, in the second terms of \eq{KA} and \eq{KB},
we have the same thing up to such replacements as $b\leftrightarrow c$ and
$i\leftrightarrow j$.
Inserting \eq{SU} into \eq{KB}, and after some algebra, we obtain its
compact and elegant representation \eq{ZY} in terms of the $_3F_2$
functions, which result from summations over the remaining indices.
\end{proof}

\noindent
{\bf Remark 4.}
Proceeding in the same fashion in the case of a $z$-dependent
$_3F_2$ function \eq{KA} and using the following
linear transformation of the Gauss hypergeometric function (see, e.g.,
\cite[p. 33, Eq. (19)]{SriMan}):
\begin{equation}\label{Lin}
_2F_1(a,b;c;z)=(1-z)^{-a}\,_2F_1\left(a,c-b;c;\frac{z}{z-1}\right)
\qquad \big(c\ne 0,-1,-2,\cdots;\; |\arg(1-z)|<\pi\big),
\end{equation}
we could simplify the second summation there
and this procedure could also lead to some interesting and useful results.

\vskip 3mm

\noindent
{\bf Remark 5.}
Since the parameters of $_2F_1$ functions in \eq{KA}
are of the form $a,\ b,\ b+1$, they are just hypergeometric representations,
via \cite[Entry (7.3.1.28)]{PBM3}, of the incomplete Beta function $B_z(a,b)$ given by
(see, e.g., \cite[Ch. 1]{SriMan})
$$
B_z(a,b)=\int_0^z\,t^{a-1}(1-t)^{b-1}\;dt.
$$
Indeed we have
$$
B_z(a,b)=a^{-1}z^a\,_2F_1(a,1-b,a+1;z)=
a^{-1}z^a(1-z)^{b-1}\,_2F_1\left(1-b,1,a+1;\frac{z}{z-1}\right),
$$
where the first equality is given by \cite[p. 263, Entry (6.6.8)]{AS} and
the second one is a result of the linear transformation \eq{Lin}.

\vskip 3mm

\noindent
{\bf Remark 6.}
Our summation formula \eq{ZY} gives an equivalent alternative
representation of \eq{KB} by Karlsson \cite{K74}.
Both \eq{ZY} and its special case \eq{ZX} are direct generalizations
of several results by other authors quoted in Section \ref{BR}.
Their equations \eq{OM}, \eq{MP}, \eq{FA} and \eq{FB}
are \emph{immediate} consequences of \eq{ZY}.

\vskip 3mm

In concluding this section, we note that, in a private communication,
Christian Krattenthaler suggested an alternative derivation of the
summation formula \eq{ZY}, which employs the equation (3.3.3)
in its limit case when $q \to 1$ along with the equation (3.1.1) of the book \cite{GaR}.

\section{Concluding remarks}

In our present investigation, we have proposed
the summation formulas \eq{ZY} and \eq{ZX} for the
Clausenian hypergeometric function $_3F_2$ with unit argument
and arbitrary negative integral parameter differences in two
pairs of the upper and lower parameters.
Our formulas are alternative representations
of two special cases of more general reduction formulas
derived in 1974 by Karlsson \cite{K74} in a form
of multiple finite sums.
The manifestly symmetric three-term relations for
$_3F_2$ functions with unit argument recorded in
\eq{ZY} and \eq{ZX} are evidently more advantageous than
some double-sum representations quoted in \eq{KB}.
They can be easily transformed, by using standard
relations for $_3F_2$ functions, according to
the specific needs of certain calculations.

As discussed already in Section \ref{sec1}, the present
summation formulas are of interest from the point of view of
practical applications in field theories in $d$-dimensional spaces
$\mathbb R^d=\mathbb R^D\oplus\mathbb R^m$, where a global rotational
$O(d)$ symmetry no more exists (see, for details, \cite{DS00,SD01,SPD05,RDS11,SP12,SP14})
or in relativistic field theories with broken Lorentz invariance
of the space-time, the so-called Lorentz violating theories.
A short review of both statistical
and high-energy physics realizations of such theories can be found in
the introduction of \cite{SP12}.

Some implications from the purely mathematical side are
mentioned in our Remarks 1 to 6.

Of course, it would be of great interest to derive
representations of similar kind for the $z$-dependent functions $_3F_2$
starting from \eq{KA} or by some other means.
A specimen of such relations is given in
\cite[Entry (7.4.1.5)]{PBM3}, it is obtained by trivially applying the
result \eq{KA} for the case when $m=n=0$:
$$
\TF{a}{b}{c}{b+1}{c+1}{z}=\frac1{c-b}
\lr{c\;\FG{a}{b}{b+1}{z}-b\;\FG{a}{c}{c+1}{z}}.
$$

Also, it should be interesting to write down relations
analogous to \eq{ZY} involving both positive and negative integers
$m$ and $n$ and to derive counterparts of equations \eq{ZY} and \eq{ZX}
for {\it bilateral} hypergeometric series \cite{Bail36}.

Finally, as suggested by the anonymous referee, similar problems might be considered
for certain special cases of Kamp{\'e} de F{\'e}riet functions,
which would be related to earlier works of Srivastava \cite{Sri78}
and Karlsson \cite{K82}.

\section*{Acknowledgements}
Our grateful thanks are due to Prof. H. W. Diehl for giving the first-named author
the possibility to work on
this paper at the Fakult\"at f\"ur Physik of the Universit\"at Duisburg-Essen
and for his warm hospitality and financial support during his stay there in 2014.
We are thankful also to Prof. T. K. Pog\'any
for reading the paper and pointing out a misprint in the first draft, and also for his numerous
friendly suggestions. Correspondence with Prof. R. B. Paris and
Prof. C. Krattenthaler is gratefully acknowledged.


\end{document}